\documentclass[12pt,reqno]{amsart}
\usepackage{colortbl} 
 \usepackage{xcolor}
\usepackage{graphics,cancel,graphicx}
\usepackage{epsfig,psfrag,manfnt}
\usepackage{bbm,subfigure,rotating,bm,float,mathdots,wasysym,scalerel}

\usepackage{tikz} 
\usetikzlibrary{arrows.meta}

\usepackage[all,cmtip]{xy}
\usepackage{amssymb,amsmath,amsfonts,amsthm,color}
\usepackage{mathrsfs}

\usepackage[new]{old-arrows}
\usepackage{stackengine}
\stackMath
\newcommand\reallywidehat[1]{%
\savestack{\tmpbox}{\stretchto{%
  \scaleto{%
    \scalerel*[\widthof{\ensuremath{#1}}]{\kern-.6pt\bigwedge\kern-.6pt}%
    {\rule[-\textheight/2]{1ex}{\textheight}}
  }{\textheight}%
}{0.5ex}}%
\stackon[1pt]{#1}{\tmpbox}%
}
\newcommand{\calA}{\mathcal{A}}

\newcommand{\mA}{\mathbb{A}}

\newcommand{\mC}{\mathbb{C}}
\newcommand{\mD}{\mathbb{D}}

\newcommand{\mF}{\mathbb{F}}

\newcommand{\mR}{\mathbb{R}}
\newcommand{\mS}{\mathbb{S}}
\newcommand{\mT}{\mathbb{T}}

\newcommand{\mZ}{\mathbb{Z}}

\newtheorem{theorem}{Theorem}[section]
\newtheorem{lemma}[theorem]{Lemma}

\newtheorem{proposition}[theorem]{Proposition}

\theoremstyle{definition}

\theoremstyle{definition}

\theoremstyle{definition}

\theoremstyle{definition}

\begin{document}

\keywords{Robust control theory, Hardy algebra, $\nu$-metric, stabilisation problem}

\subjclass[2010]{Primary 93D21; Secondary 93D15, 93B36, 30H05}

\title[]{A new $\nu$-metric computational example for the diffusion equation with \\
  boundary control and point observation}

\author[]{Amol Sasane}
\address{Department of Mathematics \\London School of Economics\\
     Houghton Street\\ London WC2A 2AE\\ United Kingdom}
\email{A.J.Sasane@LSE.ac.uk}
 
\maketitle

 \vspace{-0.9cm}
 
\begin{abstract} 
The $\nu$-metric used in robust control is computed for  control
systems with parametric uncertainty, governed by a diffusion equation in a bounded  one-dimensional 
spatial region with boundary control and point observation.
\end{abstract}

\section{Introduction}

\noindent The aim of this article is to present an example of the
computation of the $\nu$-metric, used in robust control, for  control
systems with parametric uncertainty governed by a diffusion equation with boundary control and
point observation.

In the  factorisation approach to linear control theory, 
one starts with a complex normed algebra $R$, meant to be 
the set of all transfer functions of stable  
linear systems, and then the field of fractions $\mF(R)$ over $R$ is thought of as the set of transfer functions of linear 
systems which are not necessarily stable.  Based on algebraic 
factorisations of the `plant' (i.e., to-be-controlled) transfer function, control theoretic problems can 
then be posed and solved (see, e.g. \cite{Vid}). For example, in
 the stabilisation problem, given $p\in\mF(R)$ (an unstable plant), 
one aims to find 
$c\in \mF(R)$ (a stabilising controller) such that 
the interconnection of $p$ and $c$ is stable, that is,  
the closed loop transfer function is stable:
 $$
 \textstyle 
  H(p,c):=\left[\begin{array}{cc} 
  {\scaleobj{0.81}{
\cfrac{\raisebox{-0.21em}{\text{$-p  c$}}}{\raisebox{0.21em}{\text{$1-pc$}}}}} & 
{\scaleobj{0.81}{\cfrac{\raisebox{-0.21em}{\text{$p$}}}{\raisebox{0.21em}{\text{$1-p c$}}}}}\\[0.21cm]
{\scaleobj{0.81}{\cfrac{\raisebox{-0.21em}{\text{$-c$}}}{\raisebox{0.21em}{\text{$1-pc$}}}}} & 
{\scaleobj{0.81}{\cfrac{\raisebox{-0.21em}{\text{$1$}}}{\raisebox{0.21em}{\text{$1-pc$}}}}}
      \end{array}\right]\in R^{2\times 2}.
$$
This problem can be solved if $p$ has a coprime factorisation, and for some cases of $R$ (e.g.,  the Hardy algebra $H^\infty$ of all complex-valued bounded and holomorphic functions in $\mC_+=\{s\in \mC:\text{Re}\;\!s>0\}$), the condition of existence of a coprime factorisation is also necessary.

\begin{center}
\begin{tikzpicture}[>={To[length=1.6mm,width=1.5mm]}, line width=0.5pt]
  \draw (0,0) circle (0.1);
  \draw (-0.065,0)--(0.065,0) (0,-0.065)--(0,0.065);
  \draw (3.22,-1.73) circle (0.1);
  \draw (3.155,-1.73)--(3.285,-1.73) (3.22,-1.795)--(3.22,-1.665);
  \draw (1.1,-0.38) rectangle (2.1,0.38) node[midway] {$p$};
  \draw[dash pattern=on 1.5pt off 1.2pt, draw={rgb,255:red,0;green,143;blue,143}] (0.942,-0.551) rectangle (1.942,0.209);
  \draw[dash pattern=on 1.5pt off 1.2pt, draw={rgb,255:red,0;green,143;blue,143}] (1.28,-0.209) rectangle (2.28,0.551);
  \draw (1.1,-2.11) rectangle (2.1,-1.35) node[midway] {$c$};
  \draw[->] (-0.79,0) -- (-0.1,0);
  \draw[->] (0.1,0) -- (1.1,0);
  \draw[->] (2.1,0) -- (3.22,0) -- (3.22,-1.63);
  \draw[->] (4.0,-1.73) -- (3.32,-1.73);
  \draw[->] (3.12,-1.73) -- (2.1,-1.73);
  \draw[->] (1.1,-1.73) -- (0,-1.73) -- (0,-0.1);
\end{tikzpicture}

  {\small Robust stabilisation of an uncertain plant.}
   \label{interconnection}
\end{center}

\vspace{-0.3cm}

\noindent In reality, the plant transfer function is computed from the differential equation model of the 
situation at hand, which in turn is obtained from a modelling procedure involving idealisations, 
simplifying assumptions, approximations, and measurement of parameters. So  
the plant transfer function is not known precisely, but serves only as an approximation of reality, and so  
knowing just a `nominal' plant transfer function $p$, the reality may be viewed as a perturbation $\widetilde{p}$ of this nominal transfer function $p$. 
One  then wishes that the stabilising controller $c$ found for the nominal $p$ stabilises not only $p$, but also
 all plants $\widetilde{p}$ `close enough' to $p$, 
and we would also like to be able to compute the radius of this 
neighbourhood (so that one knows the tolerated parametric uncertainty). 
  The question of what is an appropriate notion of distance one should use to measure 
closeness of unstable plants thus arises naturally. Based on the goals described  above, 
 the metric $d$ on the set of stabilisable plants should be such that  
 $d$ is easily computable, and such that 
   stabilisability is 
 a robust property of stabilisable plants in the metric $d$, i.e., if $p$ is stabilised by a controller $c$, 
 then there exists an $r>0$ such that every $\widetilde{p}$ in the ball $B(p,r):=\{\widetilde{p}:d(\widetilde{p},p)<r\}$ is stabilised by $c$. 
 Such a metric, called the $\nu$-metric, was introduced by Vinnicombe in \cite{Vin} when $R=H^\infty\cap \mC(s)$ (rational functions without poles 
in $\mC_+\cup i\mR$), and the definition of the $\nu$-metric 
 was given in terms of normalised coprime factorisations. 
 
 An abstract extension of Vinnicombe's metric was given in \cite{BalSas}  
 in order to cover the case when $R\neq H^\infty\cap \mC(s)$; rather $R$ is the class of stable transfer functions 
 of infinite dimensional systems, e.g. when $R$ is the Callier-Desoer algebra $\calA_+$ (complex Borel measures on $\mR$ with support in $[0,\infty)$ without a singular nonatomic part, with pointwise operations, multiplication given by convolution, and the total variation norm taken as the Banach algebraic norm). The starting point in \cite{BalSas} 
 was abstract, where it was assumed that $R$  satisfies some mild assumptions, and then an abstract $\nu$-metric was defined on the subset of $\mF(R)$ consisting of all elements  possessing  
 a normalised coprime factorisation. Specific examples of $R$ satisfying the abstract assumptions and the resulting particular examples 
 of the $\nu$-metric were given in \cite{BalSas}, and while the case  $R\!=\!\calA_+$ was covered, the case $R\!=\!H^\infty(\mD)$ (set of all bounded and holomorphic functions in $\mD\!=\!\{z\!\in\! \mC\!:\!|z|\!<\!1\}$) was not covered there. The case $R\!=\!H^\infty(\mD)$ was covered later  in \cite{Sas0} and \cite{FreSas}. 
 However, in   \cite{BalSas}, the $\nu$-metric relied on the existence of {\em normalised} coprime factorisations, which was a limitation, because firstly, the existence of a coprime factorisation of a plant does not 
 automatically imply the existence of a normalised one, and secondly, even if a normalised one is known to exist, there may not 
 be a computational procedure to find it \cite{PS}. On the other hand, the problem of finding coprime factorisations is much more tractable 
 than finding normalised coprime ones. In \cite{Sas},  an extension of the abstract $\nu$-metric from \cite{BalSas} was given in the single-input single-output case, circumventing the need to work
with normalised coprime factorisations.

\noindent {\bf Motivation.} The reason behind using a metric whose definition is as involved as the one recalled in Section~\ref{Section_2} is that the plants for which feedback is needed are typically {\em unstable}, 
and for such plants the 
most obvious candidate, $\|p-\widetilde{p}\|_\infty$, is not even defined. Indeed, for the plant considered here, 
we have that $p_a(s)\!\sim\! \frac{1}{s}$ as $s\to 0$, and so $p_a\not\in H^\infty$. Moreover, closeness in $\|\cdot\|_\infty$ is not the appropriate notion for feedback, whereas closeness in the $\nu$-metric $d$ is:
it was shown in \cite[Theorem~4.5]{Sas} (just as in \cite{Vin} for the rational case) that if a controller $c$ stabilises $p$, then it also stabilises every $\widetilde{p}$ with $d(p,\widetilde{p})\!<\!\mu_{p,c}$, 
where $\mu_{p,c}$ is the stability margin of the loop. 
Recall that, for $p,c\in\mF(H^\infty)$,  $\mu_{p,c}\!:=\!\frac{1}{\|H(p,c)\|_\infty}$ if $c$ stabilises $p$, and $\mu_{p,c}\!:=\!0$ otherwise. Here $\|H(p,c)\|_\infty$ is the supremum over $s\in \mC_+$ of the operator norm of  $(H(p,c))(s):\mC^2\to\mC^2$, when $\mC^2$ is endowed with the usual Euclidean norm. 
Thus a controller designed for the nominal plant $p$ with a large margin tolerates all perturbations $\widetilde{p}$ in the ball of radius $\mu_{p,c}$; this is how the metric is used in robust controller design 
(one designs $c$ to maximise $\mu_{p,c}$, and then reads off the tolerated uncertainty as a ball in the metric $d$).  
Also, the definition of $d$ in \cite{Sas} uses only (not necessarily normalised) coprime factorisations in the single-input single-output case. These are explicit 
and elementary for the example in this article, while normalised ones are generally not available (see, e.g., \cite{PS}). 
 To be of practical use, the metric must be computable for concrete plants, and we carry this out here for a system governed by a PDE. The uncertain 
parameter is the sensor location $a$ (e.g., the position of a thermocouple in a rod), and $d(p_a,p_{\widetilde{a}})$ quantifies how large an error in the sensor position 
a stabilising controller can be expected to tolerate. 
The relation between the transfer function $p_a$ and the control problem for the PDE is explained in Section~\ref{Section_3}.

\noindent {\bf Contributions.} Besides computing the $\nu$-metric for a system governed by a PDE, the article contains: 
\begin{itemize}
\item[${\scaleobj{0.72}{\bullet}}$] explicit moduli of continuity (Lipschitz estimates) 
for the coprime factors $n_a,d_a$ (Propositions~\ref{prop_raymond} and \ref{prop_sasane}); 
\item[${\scaleobj{0.72}{\bullet}}$] a verification that the index condition holds for $\widetilde{a}$ close to $a$, with an explicit neighbourhood (Lemma~\ref{lem_re} and Theorem~\ref{thm_quant}); 
\item[${\scaleobj{0.72}{\bullet}}$] a  bound of $d(p_a,p_{\widetilde{a}})$ in terms of $|a-\widetilde{a}|$ (Theorem~\ref{thm_quant}); 
and 
\item[${\scaleobj{0.72}{\bullet}}$] a closed form for the function whose supremum gives $d$, together with a numerical evaluation of this supremum (Section~\ref{Section_4}). 
\end{itemize}
We list some examples from the earlier works of the computation of the $\nu$-metric $d$ from \cite{Sas}, in the
context of delay-differential equations.  

In
\cite{Sas0}, it was shown that if $a,\widetilde{a}$ are real numbers such that $|a-\widetilde{a}|$ is small
enough, and $T>0$, then
$$
\textstyle d(e^{-sT}{\scaleobj{0.81}{
\cfrac{\raisebox{-0.21em}{\text{$s$}}}{\raisebox{0.21em}{\text{$s-a$}}}}},
e^{-sT}{\scaleobj{0.81}{
\cfrac{\raisebox{-0.21em}{\text{$s$}}}{\raisebox{0.21em}{\text{$s-\widetilde{a}$}}}}})
={\scaleobj{0.81}{
\cfrac{|a-\widetilde{a}|}{\sqrt{2}(a+\widetilde{a})}}},
$$
 quantifying the effect of uncertainty in the {\em pole}
 location. Also, it was shown that if $T\neq \widetilde{T}$ are
 distinct positive numbers, then
 $$
\textstyle d(e^{-sT}{\scaleobj{0.81}{
\cfrac{\raisebox{-0.21em}{\text{$s$}}}{\raisebox{0.21em}{\text{$s-a$}}}}},
e^{-s\widetilde{T}}{\scaleobj{0.81}{
\cfrac{\raisebox{-0.21em}{\text{$s$}}}{\raisebox{0.21em}{\text{$s-a$}}}}})=1,
$$
quantifying the effect of the uncertainty in the {\em
  delay}. In \cite{FreSas}, it was shown that if $|b|$ is small enough
then
$$
\textstyle d(e^{-sT}{\scaleobj{0.81}{
\cfrac{\raisebox{-0.21em}{\text{$s-0$}}}{\raisebox{0.21em}{\text{$s-a$}}}}},
e^{-sT}{\scaleobj{0.81}{
\cfrac{\raisebox{-0.21em}{\text{$s-b$}}}{\raisebox{0.21em}{\text{$s-a$}}}}})={\scaleobj{0.81}{
\cfrac{|b|}{\sqrt{b^2+a^2}}}},
$$
 quantifying the effect of the uncertainty in the {\em zero}
location. 

Finally, consider the transfer function
 $$
\textstyle p_\delta:={\scaleobj{0.81}{\cfrac{1}{s-(1+\delta)e^{-s}}}}
$$ 
 associated with the retarded
delay differential equation
$$
\textstyle 
\begin{array}{rcl}
{\scaleobj{0.81}{
\cfrac{\raisebox{-0.21em}{\text{$dx$}}}{\raisebox{0.21em}{\text{$dt$}}}}}(t)\!\!\!&=&\!\!\!(1+\delta) x(t-1)+u(t)\\[0.1cm]
y(t)\!\!\!&=&\!\!\! x(t)\end{array}\bigg\}, \quad t\geq 0.
$$
 Let the nominal transfer function correspond to $\delta\!=\!0$.  In
\cite{PS}, it was shown that an algebraic spectral factorisation
leading to normalised coprime factors is not possible. For $\delta$
small enough, it was shown in \cite{Sas} that
$$
\textstyle
d(p_0, p_\delta)={\scaleobj{0.81}{
\cfrac{|\delta|}{\sqrt{2(1+(1+\delta)^2)}}}}.
$$

\noindent In the present article we compute the $\nu$-metric for control systems with parametric uncertainty, when the 
control system is described by a  partial differential equation, namely  the diffusion
equation in one spatial dimension, with boundary control and point observation. 
The organisation of the article is as follows:
 \begin{itemize}
\item[${\scaleobj{0.72}{\bullet}}$]
  In Section~\ref{Section_2},
we recall the $\nu$-metric from \cite{Sas} when $R=H^\infty(\mD)$.

 \item[${\scaleobj{0.72}{\bullet}}$]
 In Section~\ref{Section_3}, we describe 
the control system  given by 
the diffusion equation in one spatial dimension with boundary control and point 
 observation.

\item[${\scaleobj{0.72}{\bullet}}$] Finally, in
Section~\ref{Section_4}, we give the details of the computation of the
$\nu$-metric for
such control systems when there is parametric uncertainty. We also give explicit moduli of continuity of the coprime factors, and a  bound for $d(p_a,p_{\widetilde{a}})$ in terms of $|a-\widetilde{a}|$.
\end{itemize}
 
 \vspace{-0.3cm}
 
\section{The $\nu$ metric: $H^\infty(\mD)$ case}
\label{Section_2}

 \vspace{-0.21cm}

\noindent We first recall the $\nu$-metric $d$ from \cite{Sas}. This section is a recap, without proofs, of material from \cite{BalSas,Sasb,Sas,FreSas}, included to make the article self-contained. What is used later is the formula for the chordal distance $\kappa$, the definition of $d$, and the properties (I3) and (I5) of the index $\iota$, which are used in Lemma~\ref{lem_re}  (to verify that the index vanishes in Section~\ref{Section_4}). 
Rather than recalling the full abstract set-up of the $\nu$-metric, we will just specialise to the case when the algebra of
stable transfer functions is given by the Hardy algebra
$H^\infty(\mD)$, that is, the Banach algebra of all bounded and
holomorphic functions defined in $\mD=\{z\in \mC:|z|<1\}$, with
pointwise operations, and the supremum norm $\|\cdot\|_\infty$
(given by $ \textstyle \|f\|_\infty:=\sup_{z\in \mD}|f(z)| $ for all
$f\in H^\infty(\mD)$).  For background on $C^*$ algebras, we refer the reader to \cite{Mur}. 

For  $r\!\in \!(0,1)$, let $\mA_r\!:=\!\{z\!\in\!
\mC\!:\!r\!<\!|z|\!<\!1\}$, and  $\text{C}_{\text{b}}(\mA_r)$ be the $C^*$-algebra consisting of
all bounded and continuous $f\!:\!\mA_r \!\to\! \mC$, equipped with
pointwise operations, the supremum norm (given by $ \textstyle \|f\|_\infty\!:=\!\sup_{z\in \mA_r}|f(z)| $ for all
$f\in \text{C}_{\text{b}}(\mA_r)$), and involution $\cdot^{*}$ defined via
pointwise complex conjugation. For $0 \!<\! r\!\le\! R \!<\! 1$, define $\pi_r^R:\text{C}_{\text{b}}(\mA_r)\to \text{C}_{\text{b}}(\mA_R)$ by
restriction: $\pi_r^R(f)\!=\!f|_{\mA_R}$ for $f\in \text{C}_{\text{b}}(\mA_r)$.
Consider the family $(\text{C}_{\text{b}}(\mA_r),\pi_r^R)$, $0\!<\!r\!\le \!R\!<\!1$. We
have $\pi_r^r$ is the identity map on $\text{C}_{\text{b}}(\mA_r)$, and
$\pi_r^R\circ \pi_\rho^r\!=\!\pi_\rho^R$ for all $0\!<\!\rho\!\le\! r\!\le\! R\!<\!1$. Now
consider the $*$-algebra $\prod_{r\in (0,1)} \text{C}_{\text{b}}(\mA_r)$, and
denote by $\mathscr{A}$ its $*$-subalgebra consisting of all 
$f\!=\!(f_r)_{r\in (0,1)}$ such that there is an index $r_0$ such that 
$\pi_r^R(f_r)\!=\!f_R$ for all $0\!<\!r_0\!\le\! r\!\le\! R\!<\!1$. Since every $\pi_r^R$
is norm decreasing, the net $(\|f_r\|_\infty)_{r\in (0,1)}$ is convergent, and we
define $\|f\|\!=\!\lim_{r\to 1} \|f_r\|_\infty$. This defines a seminorm
on $\mathscr{A}$ that satisfies the $C^*$-norm identity $\|f^* f\|\!=\!\|f\|^2$,
where $\cdot^*$ is the involution (obtained from complex conjugation). 
If $\mathscr{N}$ is the kernel of $\|\cdot\|$, then the quotient $\mathscr{A}/\mathscr{N}$ is a
$C^*$-algebra (and the norm is again denoted by $\|\cdot\|$), called the {\em direct/inductive limit of
  $(\text{C}_{\text{b}}(\mA_r),\pi_r^R)$}, and is denoted by $\varinjlim
\text{C}_{\text{b}}(\mA_r)$. To every element $f \in
\text{C}_{\text{b}}(\mA_{r_0})$, one associates a sequence $f_1\!=\!(f_r)_{r\in (0,1)}$ in
$\mathscr{A}$, where $f_{r}\!=\! 0$ if $0\!<\!r\!<\!r_{0}$, and $\pi_{r_{0}}^{r}(f)$ if
$r_{0}\!\leq\! r\!<\!1.  $ We also define
$\pi_{r}:\text{C}_{\text{b}}(\mathbb{A}_{r})\rightarrow
\varinjlim\text{C}_{\text{b}}(\mathbb{A}_{r})$ by $
\pi_{r}(f)\!:=\![f_{1}]$ for all $ f\in \text{C}_{\text{b}}(\mathbb{A}_{r}), $
where $[f_{1}]$ denotes the equivalence class in $\varinjlim
\text{C}_{\text{b}}(\mathbb{A}_{r})$ which contains $f_{1}$. Then $\pi_{r}$ are $\ast$-homomorphisms. These
maps are compatible with the connecting maps $\pi_{r}^{R}$ in the
sense that every diagram shown below is commutative.
\vspace{-0.21cm}
$$
\textstyle 
\xymatrix{
\text{C}_{\text{b}}(\mA_r) \ar[r]^{\pi_r^R} \ar[rd]_{\pi_r} & \text{C}_{\text{b}}(\mA_R) \ar[d]^{\pi_R} \\ &
 {\varinjlim \text{C}_{\text{b}}(\mA_r)}
}
$$
Then $\varinjlim \text{C}_{\text{b}}(\mathbb{A}_{r})$ is a $C^\ast$-algebra. The multiplicative identity arises from
the constant function $f\equiv1 $ in $\text{C}_{\text{b}}(\mathbb{A}_{0})$,
i.e., $\pi_{0}(f)$. The involution in  $\varinjlim 
\text{C}_{\text{b}}(\mathbb{A}_{r})$ is induced from the involution in $\text{C}_{\text{b}}(\mathbb{A}_{r})$.  

The map $H^\infty(\mD)\owns f\mapsto\pi_{0}(f)\in \varinjlim \text{C}_{\text{b}}
(\mathbb{A}_{r})$ gives an embedding of the Hardy algebra 
$H^{\infty}(\mD)$ in $\varinjlim \text{C}_{\text{b}}(\mathbb{A}_{r})$.  

For a compact Hausdorff topological space $X$, let $\text{C}(X)$ denote the $C^*$-algebra consisting of all 
continuous $f:X\to \mC$, with pointwise defined operations and involution, and the supremum norm
 (given by $\|f\|_\infty:=\sup_{x\in X}|f(x)|$ for all $f\in \text{C}(X)$). For a $C^*$-algebra $S$, 
$\text{inv}\;\! S$  denotes the set of all invertible elements of $S$.

Let $\mT\!=\!\{z\in \mC:|z|\!=\!1\}$.
For $f\!\in\!{\text{inv}}\;\!
\text{C}_{\text{b}}(\mathbb{A}_{\rho})$ and $0\!<\!\rho\!<\!r\!<\!1$, 
 defining
$f^{r}\in \text{inv}\;\!\text{C}(\mathbb{T})$ by $f^{r}(\zeta)\!=\!f(r\zeta)$
for all $\zeta\in\mathbb{T}$, we have that $f^{r}$ has a
winding number $w(f^r)$ with respect to $0$.  Set
$w(f)\!=\!w(f^{r})\!\in\!\mathbb{Z}$ for $ f\!\in\!{\text{inv}}\;\!
\text{C}_{\text{b}}(\mathbb{A}_{\rho})$. Then $w$ is well-defined.  Define
$\iota\!:\!{\text{inv}}\;\! (\varinjlim
\text{C}_{\text{b}}(\mathbb{A}_{r}))\!\rightarrow\! \mathbb{Z}$ by 
$$
\textstyle
 \iota(f)=\lim\limits_{r\rightarrow1}w(f_{r})\text{ for all  }
f=[(f_{r})]\in{\text{inv}}\;\! (\varinjlim \text{C}_{\text{b}}(\mathbb{A}
_{r})).
$$
Then $\iota$ is well-defined and satisfies the following:

\noindent (I1) $\iota(ab)= \iota (a) +\iota(b)$ for all $a,b \in
\text{inv}\;\! (\varinjlim \text{C}_{\text{b}}(\mathbb{A}_{r}))$.

\noindent(I2) $\iota(a^*)=-\iota(a)$ for all $a\in \text{inv}\;\!
(\varinjlim \text{C}_{\text{b}}(\mathbb{A}_{r}))$.

\noindent(I3) $\iota$ is locally constant.

\noindent(I4) $f\!\in\! H^\infty(\mD) \cap \text{inv} (\varinjlim
\text{C}_{\text{b}}(\mathbb{A}_{r}))$ is in $\text{inv} \;\!H^\infty(\mD)$ if
and only if $\iota(f)\!=\!0$.

\noindent(I5) If $f\in \varinjlim \text{C}_{\text{b}}(\mathbb{A}_{r})$ and
$f>0$ (in the $C^*$ algebra), then $\iota(f)=0$.

\noindent 
It was also shown in \cite{FreSas} that $\varinjlim
    \text{C}_{\text{b}}(\mathbb{A}_{r}) $ is isometrically isomorphic to $
\text{C}(\beta\mA_0 \!\setminus\! \mA_0)$ (where $
\beta \mA_0$ denotes the
Stone-\v{C}ech compactification of
$\mA_0$, which is  the maximal ideal space of the Banach algebra
$\text{C}_{\text{b}}(\mA_0)$ of all complex-valued bounded
continuous functions on $\mA_0$), and  
$\varinjlim \text{C}_{\text{b}}(\mathbb{A}_{r}) $ is a
sub-$C^\ast$-algebra of $\text{L}^\infty(\mT)$ (essentially bounded
complex-valued functions on $\mT$ with pointwise operations and the essential supremum norm).

Let $\mF(H^\infty(\mD))$ be the field of fractions of $H^\infty(\mD)$.
For $p \in \mF(H^\infty(\mD))$, a factorisation 
$$
p={\scaleobj{0.81}{
\cfrac{\raisebox{-0.21em}{\text{$n$}}}{\raisebox{0.21em}{\text{$d$}}}}},
$$
 where
$n \in H^\infty(\mD)$, $d\in H^\infty(\mD)\setminus\{ 0\}$, is a {\em
  coprime factorisation} of $p$ if there exist $x, y\in H^\infty(\mD)$
such that 
 $
n x+d y=1.
$ 
  Let $\mS(H^\infty(\mD))$ denote the set of
elements in $\mF(H^\infty(\mD))$ possessing a coprime factorisation.

The maximal ideal space of a commutative $C^*$-algebra $S$ is denoted
by $M(S)$. For $x\in S$, let $\widehat{x}$ denote the Gelfand
transform of $x$. By the Gelfand-Naimark Theorem, 
$$
\textstyle
 \|x\|_S =
\max\limits_{\varphi \in M(S)}|\widehat{x}(\varphi)|. 
$$
Let $p_1,p_2 \in \mS(H^\infty(\mD))$ have the coprime factorisations
$$
\textstyle
p_1={\scaleobj{0.81}{\cfrac{\raisebox{-0.21em}{\text{$n_1$}}}{\raisebox{0.21em}{\text{$d_1$}}}}}\quad \text{ and }\quad 
p_2= {\scaleobj{0.81}{\cfrac{\raisebox{-0.21em}{\text{$n_2$}}}{\raisebox{0.21em}{\text{$d_2$}}}}}.
$$
Then the {\em chordal distance} $\kappa(p_1,p_2)$ is
$$ 
\textstyle
\kappa(p_1,p_2)
:=
{\scaleobj{0.81}{
\sup\limits _{\varphi \in M( \varinjlim C_{\text{b}}(\mathbb{A}_{r}) )} 
\cfrac{|\widehat{n_1}(\varphi)\widehat{d_2}(\varphi)
-\widehat{n_2}(\varphi)\widehat{d_1}(\varphi)|}{
\sqrt{|\widehat{n_1}(\varphi)|^2+|\widehat{d_1}(\varphi)|^2} 
\sqrt{|\widehat{n_2}(\varphi)|^2+|\widehat{d_2}(\varphi)|^2}} .}}
$$
\noindent  
The function $\kappa$ given by the above expression is well-defined,
that is, it does not depend on the choice of coprime factorisations
for each of the plants (see \cite{Sasb}).  Moreover the above matches
the usual $\|\cdot\|_\infty$-norm of the function on the right-hand
side in $\text{L}^\infty(\mT)$ (see \cite{FreSas}).

The $\nu$ metric between $p_1,p_2$ is defined by
$$
\textstyle
d (p_1,p_2 )=\Bigg\{
\begin{array}{cl}
  \kappa(p_1,p_2) &
  \textrm{if } n_1^* n_2+d_1^* d_2 \in \text{inv}\;\!
  (\varinjlim C_{\text{b}}(\mathbb{A}_{r})) \textrm{ and }
  \\
  &
  \phantom{\text{if }}\iota ( n_1^* n_2+d_1^* d_2)=0, \\[0.1cm]
  1 & \textrm{otherwise}. \end{array}
$$
 It was shown in \cite{Sas} that $d$ is a metric on
$\mS(H^\infty(\mD))$ making stabilisability a robust property.

For our control system described by a partial differential equation,
the transfer function is obtained by taking the one-sided Laplace
transform, and so it will be natural to use the Hardy algebra
$H^\infty$ of the half-plane $\mC_+\!:=\!\{s\!\in\! \mC\!:\!\text{Re}\;\!s\!>\!0\}$,
consisting of all bounded and holomorphic functions $f$ defined on
$\mC_+$. As the index $\iota$ is more naturally expressed in $\mD$, we
presented the $\nu$-metric for $H^\infty(\mD)$ above. But one
can transplant data from the half-plane to the unit disc and vice
versa using the M\"obius transforms 
$$
\mC_+\owns
s={\scaleobj{0.81}{\cfrac{\raisebox{-0.21em}{\text{$1+z$}}}{\raisebox{0.21em}{\text{$1-z$}}}}}
\leftrightarrow {\scaleobj{0.81}{\cfrac{\raisebox{-0.21em}{\text{$s-1$}}}{\raisebox{0.21em}{\text{$s+1$}}}}}=z\in \mD,
$$ as done in
the following section in the course of the computation. 

In Section~\ref{Section_4}, we will need to check that the index $\iota$ of the element $n_1^*n_2+d_1^*d_2$ vanishes. The following lemma gives a convenient sufficient condition for this, which can be checked directly: the real part of the element must be bounded below by a positive constant.

\begin{lemma}
\label{lem_re}
Let $f\in \varinjlim \text{\em C}_{\text{b}}(\mA_r)$ be such that for some $\epsilon>0,$
$$
\textstyle \text{\em Re}\;\!f:=\frac{1}{2}(f+f^*)\ge \epsilon\cdot 1,
$$ 
that is$,$ $\text{\em Re}\;\!f-\epsilon\cdot 1$ is a positive element of the $C^*$-algebra $
\varinjlim \text{\em C}_{\text{b}}(\mA_r)$. 
Then $f\in \text{\em inv}\;\!(\varinjlim \text{\em C}_{\text{b}}(\mA_r))$ and $\iota(f)=0$. 
\end{lemma}
\begin{proof}
As $\varinjlim \text{C}_{\text{b}}(\mA_r)$ is a unital sub-$C^*$-algebra of $\text{L}^\infty(\mT)$, it is inverse closed (spectral permanence, see \cite[Theorem~2.1.11]{Mur}), and positivity in $\varinjlim \text{C}_{\text{b}}(\mA_r)$ is the same as 
positivity in $\text{L}^\infty(\mT)$. For $t\in[0,1]$, let 
$$
\textstyle f_t\!:=\!\text{Re}\;\!f+i(1-t)\,\text{Im}\;\!f,\;\text{  where }\;\text{Im}\;\!f\!:=\!\frac{1}{2i}(f-f^*).
$$ 
Then $f_t\in \varinjlim \text{C}_{\text{b}}(\mA_r)$, $f_0\!=\!f$, $f_1\!=\!\text{Re}\;\!f$, and $t\mapsto f_t$ is norm continuous. 
We have $|f_t|\!\ge\! |\text{Re}\;\!f_t|\!=\!\text{Re}\;\!f\!\ge\!\epsilon$ almost everywhere on $\mT$. Thus $f_t$ is invertible in $\text{L}^\infty(\mT)$, and hence in $\varinjlim \text{C}_{\text{b}}(\mA_r)$. 
By (I3), $t\mapsto \iota(f_t)$ is locally constant on the connected set $[0,1]$, and so is in fact constant. Since $f_1\!=\!\text{Re}\;\!f\ge\epsilon\!>\!0$, (I5) gives $\iota(f_1)\!=\!0$. So $\iota(f)\!=\!\iota(f_0)\!=\!0$.
\end{proof}

\noindent The positive lower bound $\epsilon$ is not superfluous: the pointwise positivity of $\text{Re}\;\!f$ does not suffice, since the infimum could be $0$. 

\section{Control system}
\label{Section_3}

\noindent In this section, we will describe 
the control system we consider, given by 
the diffusion equation in one spatial dimension with boundary control and point 
observation. For details on viewing this as a linear control system with an  infinite-dimensional state space, we refer the reader to  
 \cite[Examples~2.3.7, 4.3.12, 7.1.11, and Exercise~7.14]{CurZwa}.

 \noindent Let $a\!\in\! (0,1)$. Consider the following 
partial differential equation (the diffusion equation) with a boundary control input given as a
Neumann boundary condition, and a point observation:
$$
\textstyle 
\begin{array}{rll}
  \text{(PDE)} \;
  & {\scaleobj{0.81}{
\cfrac{\raisebox{-0.21em}{\text{$\partial w$}}}{\raisebox{0.21em}{\text{$\partial t$}}}}} (x,t)
  \!=\!{\scaleobj{0.81}{
\cfrac{\raisebox{-0.21em}{\text{$\partial^2 w$}}}{\raisebox{0.21em}{\text{$\partial x^2$}}}}}(x,t) 
& \;\; (0\!<\!x\!<\!1, \;\!t\!\ge\! 0)\\[-0.06cm]
  \text{(BC)} \;
  & {\scaleobj{0.81}{
\cfrac{\raisebox{-0.21em}{\text{$\partial w$}}}{\raisebox{0.21em}{\text{$\partial x$}}}}}(0,t)\!=\!0
  \text{ and } 
  {\scaleobj{0.81}{
\cfrac{\raisebox{-0.21em}{\text{$\partial w$}}}{\raisebox{0.21em}{\text{$\partial x$}}}}}(1,t)\!=\!u(t) 
& \;\;(t\!\ge\! 0)\\[0.21cm]
\text{(IC)} \;
& w(x,0)=0 
&\;\;(0\!<\!x\!<\!1)\\[0.15cm]
\text{(Observation)} \;& y(t)\!=\!w(a, t) 
&\;\; (t\!\ge\! 0).
\end{array}
$$
This control system is described by the transfer function given by
 $$
 \textstyle 
 p_a(s)=
 {\scaleobj{0.81}{
\cfrac{\cosh (a \sqrt{s} )}{\sqrt{s}\cdot \sinh \sqrt{s}}}}.
 $$ 
 \noindent For $s\in
 \mC\setminus \{0\}$,  $\sqrt{s}:=e^{\frac{1}{2}\textrm{Log}\;\!s}$, where
 $\textrm{Log}\;\!s:=\log |s|+i\textrm{Arg}\;\!s$ is the principal
 logarithm of $s$, and $\textrm{Arg}\;\!s$ is the unique number in
 $(-\pi,\pi]$ such that $s=|s|(\cos (\textrm{Arg}\;\!s)+i\sin
   (\textrm{Arg}\;\!s))$. Also, $\sqrt{0}:=0$.

\noindent {\bf Relation to the control problem.}
Taking Laplace transform with respect to the $t$ variable (noting the initial condition (IC) $w(\cdot,0)\!=\!0$), 
and denoting the transforms of $t\mapsto w(x,t)$ and $u$ by $\widehat{w}(x,s)$ and $\widehat{u}(s)$, respectively, 
 the partial differential equation (PDE), and the boundary conditions (BC) yield 
$$
\textstyle 
\begin{array}{lll} 
s\widehat{w}=\widehat{w}''\text{ (where }\cdot'\text{ denotes taking derivative with respect to }x\text{)}, \\[0.06cm]
\widehat{w}'(0,s)=0,\text{ and }\widehat{w}'(1,s)=\widehat{u}(s).
\end{array}
$$ 
The solutions of the resulting ODE with $\widehat{w}'(0,s)\!=\!0$ have the form 
 $
\textstyle \widehat{w}(x,s)=A(s)\cosh(\sqrt{s}\,x),
$ 
 and the boundary condition at $x\!=\!1$ gives 
$A(s)\sqrt{s}\sinh\sqrt{s}=\widehat{u}(s)$ (note $\sinh\sqrt{s}\neq 0$ for $s\in\mC_+$, as its zeroes are at $\sqrt{s}\in i\pi\mZ$). Hence 
$$
\widehat{y}(s)=\widehat{w}(a,s)=p_a(s)\widehat{u}(s),
$$
that is, $p_a$ is the transfer function of the input-output map $u\mapsto y$ (heat flux applied at $x\!=\!1$ $\mapsto$ temperature measured at $x\!=\!a$). 
It is holomorphic in $\mC_+$ but not bounded, since $p_a(s)\sim \frac{1}{s}$ as $s\to0$. 
 Hence $p_a\notin H^\infty$, and the plant is unstable, i.e., $p_a\in \mF(H^\infty)\setminus H^\infty$. 
 The control problem for the PDE is to find a controller $c$ stabilising $p_a$, as in the Introduction (with $R=H^\infty$). 
If $p_a=\frac{n_a}{d_a}$ is a coprime factorisation, with $n_a f+d_a g=1$ for some $f,g\in H^\infty$ and $g\neq 0$, then the controller $c=-\frac{f}{g}\in \mF(H^\infty)$ 
stabilises $p_a$. 
Thus the coprime factorisation constructed in Section~\ref{Section_4} is what is needed both to solve the stabilisation problem for the PDE, and to compute the $\nu$-metric. 
The parameter $a$ is the (possibly imperfectly known) sensor position, so that a perturbation $\widetilde{a}$ of $a$ gives a perturbed plant $p_{\widetilde{a}}$.

\section{Computation of the $\nu$ metric $d(p_a,p_{\widetilde{a}})$} 
\label{Section_4}

\vspace{-0.21cm}

\noindent In this section, we give the details of the computation of the
$\nu$-metric $d(p_a,p_{\widetilde{a}})$ between the transfer functions $p_a$ and $p_{\widetilde{a}}$ of control systems 
as in the previous section. We begin by finding coprime factorisations of $p_a$ and $p_{\widetilde{a}}$.

Define $n_a,d_a$ by 
 $$
 \textstyle 
 \begin{array}{rcl}
 n_a(s)\!\!\!&=&\!\!\!{\scaleobj{0.81}{
\cfrac{\raisebox{-0.21em}{\text{$1$}}}{ \sqrt{s}+1}}}\cdot  {\scaleobj{0.81}{
\cfrac{\sinh (a\sqrt{s})}{\sinh \sqrt{s}}}}
\quad\text{  and }\quad\\[0.3cm]
 d_a(s)\!\!\!&=&\!\!\!{\scaleobj{0.81}{
\cfrac{\sqrt{s}}{\sqrt{s}+1}}} \cdot {\scaleobj{0.81}{
\cfrac{\sinh (a\sqrt{s})}{\cosh (a\sqrt{s})}}}.
\end{array}
 $$ 
 In what follows below, it will be convenient to use the angular domain $ \Delta\!=\!\{z\!\in\! \mC\!:\! |\text{Arg}\;\!
 z|\!<\!\frac{\pi}{4}\}$, and the function $R_a$ defined by 
 $$
 R_a(z):=
 {\scaleobj{0.81}{
\cfrac{\raisebox{-0.0em}{\text{$\sinh(az)$}}}{\raisebox{0.0em}{\text{$\sinh z$}}}}}, \quad z\in \Delta.
$$

\begin{proposition}
\label{241251440}
$n_a,d_a\in H^\infty$. 
\end{proposition}
\vspace{-0.39cm}
\begin{proof} It suffices to make the substitution 
$z\!:=\!\sqrt{s}$ where $z\in \Delta$. We note that the meromorphic function
 $$
 \textstyle 
 R_a(z)=
{\scaleobj{0.81}{
\cfrac{\sinh(az)}{\sinh z}}},
$$
has a removable singularity at $z\!=\!0$. It is thus enough to show the boundedness for $z\in \Delta$ with $\text{Re}\;\! z\!>\!\delta\!>\!0$, which follows by writing  
$$
\textstyle R_a(z)=e^{-(1-a)z} {\scaleobj{0.81}{
\cfrac{\raisebox{-0.21em}{\text{$1-e^{-2az}$}}}{\raisebox{0.21em}{\text{$1-e^{-2z}$}}}}},
$$ 
and noting that the numerator $|1-e^{-2az}|\!\le\! 2$, and the denominator $|1-e^{-2z}|\!\ge\! 1-e^{-2\;\!\text{Re}\;\!z}
\!\ge\! 1-e^{-2\delta}\!>\!0$. 
As $f_a(z)\!:=\!(1+z)^{-1}R_a(z)\to 0$  as $z\to \infty$ in
 $\Delta$,  $f_a$ is bounded in $\Delta$. Consequently, $n_a\in H^\infty$.

 For $d_a$, we proceed as follows. The map $\Delta\!\owns\! z\!\mapsto\! {\scaleobj{0.81}{
\cfrac{\raisebox{-0.21em}{\text{$z$}}}{\raisebox{0.21em}{\text{$z+1$}}}}}
$ is bounded as
 $$
 |z+\!1|^2\!=\!(\text{Re}\;\! z+\!1)^2\!+(\text{Im}\;\!z)^2\!=\!|z|^2\!+\!1+2\;\!\text{Re}\;\!z\!\ge\! |z|^2\!+\!1+0\!\ge\! |z|^2. 
 $$ 
 We claim $\tanh(a\;\!\cdot )$ is also bounded in $\Delta$. As $\tanh(az)\to 0$ as $z\to 0$, it is enough to show the boundedness for $z\in \Delta$ with $\text{Re}\;\! z\!>\!\delta\!>\!0$. 
 But 
 $$
 \textstyle 
 \tanh (az)={\scaleobj{0.81}{
\cfrac{\raisebox{-0.21em}{\text{$1-e^{-2az}$}}}{\raisebox{0.21em}{\text{$1+e^{-2az}$}}}}},
 $$
 and $|1\!-\!e^{-2az}|\!\le\! 2$, while  $|1\!+\!e^{-2az}|\!\ge\! 1\!-\!e^{-2a\;\!\text{Re}\;\!z}
\!\ge\! 1\!-\!e^{-2a\delta}\!>\!0$. 
 Thus 
 $$
 \textstyle
g_a(z):={\scaleobj{0.81}{
\cfrac{\raisebox{-0.21em}{\text{$z$}}}{\raisebox{0.21em}{\text{$z+1$}}}}}\tanh(az)
$$  
is bounded in $\Delta$. So $d_a\in H^\infty$.
  \end{proof}
 
  \begin{proposition}
 \label{28_1_25_1011}
 $n_a,d_a$ are coprime in $H^\infty$.
\end{proposition}
\begin{proof}  
Let $\overline{\Delta}\!=\!\{z\in\mathbb{C}: |\operatorname{Arg} z|\!\le\!\frac{\pi}{4}\}\cup\{0\}$, and 
$
\textstyle R_a(z)\!=\!\frac{\sinh(az)}{\sinh z}.
$ 
The restriction of the map $s\mapsto z\!=\!\sqrt{s}$ to 
$\overline{\mathbb{C}_+}\!:=\!\{s\in \mC:\text{Re}\;s\!\ge\! 0\}$ is a homeomorphism onto $\overline{\Delta}$, and with $z\!=\!\sqrt{s}$, we have
$$
\textstyle 
n_a(s)={\scaleobj{0.81}{
\cfrac{R_a(z)}{1+z}}},\qquad d_a(s)={\scaleobj{0.81}{
\cfrac{z}{z+1}}}\tanh(az), \quad z=\sqrt{s}.
$$
Since $\operatorname{Re} z\!\ge\! 0$ on $\overline{\Delta}$, we obtain $|1+z|\!\ge\! 1$.

\smallskip

\noindent {\em Step 1: Continuity on $\overline{\mathbb{C}_+}$ and zeroes.} The zeroes of $\sinh $ are $z\!=\!i\pi k$ and those of $\sinh(a\cdot)$ are $z\!=\!\frac{i\pi k}{a}$ ($k\in\mathbb{Z}$), 
while the zeroes of $\cosh(a\cdot)$ are 
$$
\textstyle z= {\scaleobj{0.81}{
\cfrac{i\pi(k+\frac12)}{a}}}\quad (k\in\mathbb{Z}).
$$ 
All of these, except $z\!=\!0$, have argument $\pm\frac{\pi}{2}$, and so do not lie in $\overline{\Delta}$. 
At $z\!=\!0$ the quotient $R_a$ has a removable singularity with value $a$. 
So $n_a,d_a$ are continuous on $\overline{\mathbb{C}_+}$ (and analytic in $\mathbb{C}_+$), with $n_a(0)\!=\!a$ and $d_a(0)\!=\!0$. 
Moreover, $n_a$ has no zeroes in $\overline{\mathbb{C}_+}$: $n_a(0)\!=\!a\!\neq\!0$, and for $z\in\overline{\Delta}\setminus\{0\}$ we have 
$\sinh z\!\neq\! 0$ and $\sinh(az)\!\neq\!0$, giving $R_a(z)\!\neq\!0$.

\smallskip

\noindent {\em Step 2: Behaviour at infinity.} For $z\in\overline{\Delta}$, we have 
$$
\textstyle x:=\operatorname{Re} z\ge {\scaleobj{0.81}{
\cfrac{|z|}{\sqrt{2}}}},
$$
and so $x\to\infty$ as $|z|\to\infty$. Then as $|z|\to\infty$, we have 
$$
\textstyle 
|R_a(z)|\le {\scaleobj{0.81}{
\cfrac{e^{ax}+e^{-ax}}{e^{x}-e^{-x}}}}=e^{-(1-a)x} {\scaleobj{0.81}{
\cfrac{1+e^{-2ax}}{1-e^{-2x}}}}\longrightarrow 0,
$$
because $a\!<\!1$. Thus using $|1+z|\!\ge\!1$, we get  $n_a(s)\to 0$ as $|s|\to\infty$. 
Also, as $|z|\to\infty$, we have $\frac{z}{z+1}\to 1$, and 
$$
\textstyle 
\tanh(az)={\scaleobj{0.81}{\cfrac{1-e^{-2az}}{1+e^{-2az}}}}\to 1
$$
(since $|e^{-2az}|=e^{-2a\,\text{Re}\,z}\to0$).
 So $d_a(s)\to 1$ as $|s|\to \infty$. 
Hence there is $r>0$ such that $|n_a(s)|^2+|d_a(s)|^2\!\ge\! \frac{1}{2}$ for all $s\in\overline{\mathbb{C}_+}$ with $|s|\!\ge\! r$. 

\smallskip

\noindent {\em Step 3: Uniform lower bound.} As $F\!:=\!|n_a|^2+|d_a|^2$ is continuous on the compact set 
$K\!=\!\{s\in\overline{\mathbb{C}_+}: |s|\le r\}$, and $F\!\ge\!|n_a|^2\!>\!0$ on $K$ by Step~1, we get  $m\!:=\!\min_K F\!>\!0$. By Step~2, 
$$
\textstyle 
|n_a(s)|^2+|d_a(s)|^2\ge \delta:=\min\{m,\tfrac{1}{2}\}\!>\!0\text{ for all } s\in\overline{\mathbb{C}_+}.
$$
As $(|n_a|+|d_a|)^2\!\ge\! |n_a|^2+|d_a|^2$, also $|n_a(s)|+|d_a(s)|\!\ge\!\sqrt{\delta}$ for all $s\in\mathbb{C}_+$. 
The corona theorem for $H^\infty$ now implies that $n_a,d_a$ are coprime. 
\end{proof}

\noindent Below we use the fact that for $z=x+iy\in \mC$, 
$$
\textstyle 
\begin{array}{rcl}
|\cosh z|^2\!\!\!&=&\!\!\!(\sinh x)^2+(\cos y)^2=(\cosh x)^2-(\sin y)^2,\text{ and}\\[0.1cm]
 |\sinh z|^2\!\!\!&=&\!\!\!(\sinh x)^2+(\sin y)^2.
 \end{array}
 $$
 So $|\cosh z|\!\le\!\cosh x$ and $|\sinh z|\!\ge\!|\sinh x|\!=\!\sinh |x|$. Also, for $z\in\Delta$, 
 $$
 \textstyle 
 x\!=\!\text{Re}\;\!z\!>\!{\scaleobj{0.81}{\cfrac{|z|}{\sqrt{2}}}}, \quad |1+z|^2\!=\!|z|^2+1+2x\!\ge\! |z|^2, \;\text{ and }\;\;|1+z|\!\ge\! 1.
$$

\begin{proposition}
 \label{prop_raymond}
 For all $a,\widetilde{a}\in (0,1),$  
 $
 \textstyle 
  \|n_a-n_{\widetilde{a}}\|_\infty\le 2|a-\widetilde{a}|.
  $
  
  \noindent  
   In particular$,$ $(0,1)\owns a\mapsto n_a\in H^\infty$ is Lipschitz continuous. 
 \end{proposition} 
 \vspace{-0.3cm}
 \begin{proof}
  The substitution $z\!=\!\sqrt{s}\in\Delta$ gives $n_a(s)\!=\!{\scaleobj{0.81}{\cfrac{1}{1+z}}}R_a(z)$. Since $\sinh z\!\neq\!0$ in $\Delta$ and 
 $$
 \textstyle {\scaleobj{0.81}{\cfrac{\partial\hfill }{\partial a}}}R_a(z)={\scaleobj{0.81}{\cfrac{z\cosh(az)}{\sinh z}}},
 $$
  we have 
 $$
 \textstyle 
 n_a(s)-n_{\widetilde{a}}(s)=\int_{\widetilde{a}}^{a}{\scaleobj{0.81}{\cfrac{z\cosh(\alpha z)}{(1+z)\sinh z}}}\;\!d\alpha. 
 $$ 
 So
 $$
 \textstyle 
 |n_a(s)-n_{\widetilde{a}}(s)|\le |a-\widetilde{a}|\sup\limits_{\alpha \in(0,1)}{\scaleobj{0.81}{\cfrac{|z||\cosh(\alpha z)|}{|1+z||\sinh z|}}}.
 $$
 If $z\in \Delta$,  then $|z|\!\le\!|1+z|$. If  $z\in \Delta$ moreover satisfies $|z|\!\ge \!1$, then  
 $$
 \textstyle x\!=\!\text{Re}\;\!z\!>\!{\scaleobj{0.81}{\cfrac{|z|}{\sqrt{2}}}}\!\ge \!{\scaleobj{0.81}{\cfrac{1}{\sqrt{2}}}}.
 $$
  Also, 
 $|\cosh(\alpha z)|\!\le\!\cosh(\alpha x)\!\le\! \cosh x$,  and $ |\sinh z|\!\ge\!\sinh x$.  
  So  
 $$
 \textstyle 
 \sup\limits_{\alpha \in(0,1)}{\scaleobj{0.81}{\cfrac{|z||\cosh(\alpha z)|}{|1+z||\sinh z|}}} 
 \le {\scaleobj{0.81}{\cfrac{1}{\tanh x}}} \!\le\!
 {\scaleobj{0.81}{\cfrac{1}{\tanh (1/\sqrt{2})}}}\le 2.
 $$ 
 If $|z|\!<\!1$, then $|\cosh(\alpha z)|\!\le\!\cosh 1$, 
 and $|1+z|\!\ge\! 1$. 
 Moreover, as
 $$
 \textstyle \sinh z=z+\sum\limits_{k\ge1}{\scaleobj{0.81}{\cfrac{z^{2k+1}}{(2k+1)!}}},
 $$ 
 we get 
 $$
 \textstyle 
 |(\sinh z)-z|\le\sum\limits_{k\ge1}{\scaleobj{0.81}{\cfrac{|z|^{2k+1}}{(2k+1)!}}}=(\sinh|z|)-|z|,
 $$ 
and so $|\sinh z|\!\ge\!|z|-|\sinh z-z|\!\ge\! 2|z|-\sinh|z|$. 
Since $\sinh$ is convex on $[0,1]$ with $\sinh0\!=\!0$, we have $\sinh r\!\le\! r\sinh1$ for $r\in[0,1]$. 
Consequently, $|\sinh z|\!\ge\!(2-\sinh1)|z|$
 for $|z|\!<\!1$.
  So 
   $$
   \textstyle 
   \sup\limits_{\alpha \in(0,1)}{\scaleobj{0.81}{\cfrac{|z||\cosh(\alpha z)|}{|1+z||\sinh z|}}}\le 
   {\scaleobj{0.81}{\cfrac{\cosh 1}{2-\sinh 1}}}\le 2.
   $$ 
 Thus 
  $
 \textstyle |n_a(s)-n_{\widetilde{a}}(s)|\le 2|a-\widetilde{a}|.
 $ 
 Taking the supremum over $s\in\mC_+$, the proof is completed. 
\end{proof}

\begin{proposition}
 \label{prop_sasane}
 For all $a,\widetilde{a}\in (0,1),$ 
 $$
 \textstyle  \|d_a-d_{\widetilde{a}}\|_\infty\le {\scaleobj{0.81}{\cfrac{\sqrt{2}}{\min\{a,\widetilde{a}\}}}} |a-\widetilde{a}|.
 $$
  In particular$,$ $(0,1)\owns a\mapsto d_a\in H^\infty$ is locally Lipschitz continuous. 
 \end{proposition}
 \begin{proof}
 With $z\!=\!\sqrt{s}\in\Delta$, we have $d_a(s)\!=\!{\scaleobj{0.81}{\cfrac{z}{1+z}}}\tanh(az)$. As 
 $$
 \textstyle {\scaleobj{0.81}{\cfrac{\partial\hfill }{\partial a}}}d_a={\scaleobj{0.81}{\cfrac{z^2}{(1+z)(\cosh(az))^2}}},
 $$
 we obtain (with $m\!:=\!\min\{a,\widetilde{a}\}$, and using $|1+z|\!\ge\!|z|$ for $z\in \Delta$) that 
 $$
 \textstyle 
 |d_a(s)-d_{\widetilde{a}}(s)|
 \le |a-\widetilde{a}|\sup\limits_{\alpha \ge m}
 {\scaleobj{0.81}{\cfrac{|z|^2}{|1+z||\cosh(\alpha z)|^2}}}
 \le |a-\widetilde{a}|\sup\limits _{\alpha\ge m}{\scaleobj{0.81}{\cfrac{|z|}{|\cosh(\alpha z)|^2}}}.
 $$
  Let $\alpha z\!=\!u+iv$, where $u,v\in \mathbb{R}$. Then
  $|v|\!<\!u$ as $\alpha z\in\Delta$. We claim 
  $$
  \textstyle 
  |\cosh (\alpha z)|^2=(\sinh u)^2+(\cos v)^2\ge 1.
  $$
   Indeed, if $u\!\le\!\frac{\pi}{2}$, then $(\cos v)^2\!\ge\! (\cos u)^2$, and so 
 $$
 \textstyle |\cosh (\alpha z)|^2\ge (\sinh u)^2+(\cos u)^2=\frac12(\cosh (2u)+\cos (2u))\ge 1,
 $$
  where the last inequality follows by using  
  $$
  \textstyle \cosh t\ge 1+{\scaleobj{0.81}{\cfrac{t^2}{2}}}\text{ and }\cos t\ge 1-{\scaleobj{0.81}{\cfrac{t^2}{2}}}\text{ for all }t\in \mR.
  $$
   If $u\!>\!\frac{\pi}{2}$, then $|\cosh (\alpha z)|^2\!\ge\! (\sinh u)^2\!>\!1$. 
 Also,
 $$
 \textstyle 
 |\cosh (\alpha z)|^2\ge (\sinh u)^2\ge u^2\ge \alpha^2{\scaleobj{0.81}{\cfrac{|z|^2}{2}}}.
 $$
  Hence for all $z\in \Delta$, 
 $$
 \textstyle 
 {\scaleobj{0.81}{\cfrac{|z|}{|\cosh(\alpha z)|^2}}}\le\min\Big\{|z|,{\scaleobj{0.81}{\cfrac{2}{\alpha^2|z|}}}\Big\}.
 $$
 If $|z|\!\le\!{\scaleobj{0.81}{\cfrac{\sqrt{2}}{\alpha}}}$, then 
  $
 \textstyle \min\Big\{|z|,{\scaleobj{0.81}{\cfrac{2}{\alpha^2|z|}}}\Big\}\!\le\!
  |z|\!\le\!{\scaleobj{0.81}{\cfrac{\sqrt{2}}{\alpha}}}.
 $

  \noindent If $|z|\!\ge\!{\scaleobj{0.81}{\cfrac{\sqrt{2}}{\alpha}}}$, then 
$
\textstyle 
\min\Big\{|z|,{\scaleobj{0.81}{\cfrac{2}{\alpha^2|z|}}}\Big\}
\!\le\!{\scaleobj{0.81}{\cfrac{2}{\alpha^2|z|}}}\!\le\!{\scaleobj{0.81}{\cfrac{\sqrt{2}}{\alpha}}}.
$

\noindent 
Thus 
$$
\textstyle 
  {\scaleobj{0.81}{\cfrac{|z|}{|\cosh(\alpha z)|^2}}}
  \le \min\Big\{|z|,{\scaleobj{0.81}{\cfrac{2}{\alpha^2|z|}}}\Big\} 
  \le{\scaleobj{0.81}{\cfrac{\sqrt{2}}{\alpha }}}
  \le{\scaleobj{0.81}{\cfrac{\sqrt{2}}{m}}}.
 $$
 Hence 
 $$
 \textstyle
  |d_a(s)-d_{\widetilde{a}}(s)|\le |a-\widetilde{a}|{\scaleobj{0.81}{\cfrac{\sqrt{2}}{m}}}
  = |a-\widetilde{a}|{\scaleobj{0.81}{\cfrac{\sqrt{2}}{\min\{a,\widetilde{a}\}}}}.
  $$
  Taking the supremum over $s\in\mC_+$, the proof is completed. 
\end{proof}

\noindent We remark that the dependence of the constant on $\min\{a,\widetilde{a}\}$ in Proposition~\ref{prop_sasane} cannot be removed: there is no constant $C\!>\!0$ such that $\|d_a-d_{\widetilde{a}}\|_\infty\!\le\! C|a-\widetilde{a}|$ for all $a,\widetilde{a}\in(0,1)$. Indeed, for each fixed $s\in\mC_+$ we have $\|d_a-d_{\widetilde{a}}\|_\infty\!\ge\!|d_a(s)-d_{\widetilde{a}}(s)|$, and $d_{\widetilde{a}}(s)\to0$ as $\widetilde{a}\searrow 0$, showing 
$$
\textstyle 
\liminf\limits_{\widetilde{a}\searrow 0}\|d_a-d_{\widetilde{a}}\|_\infty\ge|d_a(s)|.
$$
Taking the supremum over $s\in\mC_+$ gives 
$$
\textstyle
\liminf\limits_{\widetilde{a}\searrow 0}\|d_a-d_{\widetilde{a}}\|_\infty\ge\|d_a\|_\infty\ge\lim\limits_{x\to\infty}{\scaleobj{0.81}{\cfrac{x}{1+x}}}\tanh(ax)=1.
$$ 
On the other hand, $|a-\widetilde{a}|\to a$ as $\widetilde{a}\searrow 0$. Hence a Lipschitz estimate $\|d_a-d_{\widetilde{a}}\|_\infty\!\le\! C|a-\widetilde{a}|$ for all $\widetilde{a}\in(0,1)$ forces $C\!\ge \!\frac{1}{a}$, and since $a\in(0,1)$ is arbitrary, this proves the claim.

  From Proposition~\ref{28_1_25_1011}, there exists a $\delta_a>0$ such that 
$$
\textstyle S_a:=|n_a|^2+|d_a|^2\ge\delta_a\text{ on }\overline{\mC_+}, \text{ and in particular on }i\mR.
$$

\begin{theorem}
\label{thm_quant}
Let $a\in(0,1),$ $\delta_a$ be as above$,$ and $\Lambda_a:=\sqrt{4+\frac{8}{a^2}}$. 

\vspace{-0.21cm}

\noindent Let $\widetilde{a}\in(0,1)$ satisfy 
 $
|a-\widetilde{a}|\le\min\Big\{{\scaleobj{0.81}{\cfrac{a}{2}}},{\scaleobj{0.81}{\cfrac{\sqrt{\delta_a}}{2\Lambda_a}}}\Big\}.
$ 
Then$\;:$

\vspace{-0.39cm}

\begin{enumerate}
\item $\text{\em Re}(n_a^*n_{\widetilde{a}}+d_a^*d_{\widetilde{a}})\ge {\scaleobj{0.81}{\cfrac{\delta_a}{2}}}$  on $i\mR$. 

\noindent So $n_a^*n_{\widetilde{a}}+d_a^*d_{\widetilde{a}}\in\text{\em inv}\;\!(\varinjlim \text{\em C}_{\text{b}}(\mA_r)),$ 
and $\iota(n_a^*n_{\widetilde{a}}+d_a^*d_{\widetilde{a}})=0$. 

\vspace{0.1cm}

\item  $d(p_a,p_{\widetilde{a}})=\kappa(p_a,p_{\widetilde{a}})\le {\scaleobj{0.81}{\cfrac{\Lambda_a}{\sqrt{\delta_a}}}}|a-\widetilde{a}|$. 
\end{enumerate}
\end{theorem}
\vspace{-0.6cm}
\begin{proof}
As $|a-\widetilde{a}|\!\le\! \frac{a}{2}$, we have $\min\{a,\widetilde{a}\}\!\ge\!  \frac{a}{2}$. By Proposition~\ref{prop_raymond} and Proposition~\ref{prop_sasane}, $\|n_a-n_{\widetilde{a}}\|_\infty\!\le\! 2|a-\widetilde{a}|$ and 
$\|d_a-d_{\widetilde{a}}\|_\infty\!\le\! \frac{2\sqrt{2}}{a}|a-\widetilde{a}|$. So  we have  
 $
\textstyle
\eta:=\sqrt{|n_a-n_{\widetilde{a}}|^2+|d_a-d_{\widetilde{a}}|^2}\le \Lambda_a|a-\widetilde{a}|\le \frac{1}{2}\sqrt{\delta_a} 
$ on $i\mR$.

\noindent 
(1) On $i\mR$, we have (with $S_a:=|n_a|^2+|d_a|^2$) that 
$$
\textstyle
\begin{array}{rcl}
\text{Re}(n_a^*n_{\widetilde{a}}+d_a^*d_{\widetilde{a}})\!\!\!&=&\!\!\!S_a+\text{Re}\;\!\pmb{(}n_a^*(n_{\widetilde{a}}-n_a)+d_a^*(d_{\widetilde{a}}-d_a)\pmb{)}\\[0.1cm]
\!\!\!&\ge&\!\!\! S_a-\sqrt{S_a}\eta\quad \text{(Cauchy-Schwarz inequality)}\\[0.1cm]
\!\!\!&=&\!\!\!\sqrt{S_a}(\sqrt{S_a}-\eta)\\[0.1cm]
\!\!\!&\ge&\!\!\!  \sqrt{\delta_a}\cdot\tfrac{1}{2}\sqrt{\delta_a}=\frac{\delta_a}{2}.
\end{array}
$$
\phantom{(1) }Lemma~\ref{lem_re} (with $\epsilon\!=\!\frac{\delta_a}{2}$) gives the invertibility and $\iota=0$. 

\noindent 
(2) By (1), $d=\kappa$. As $n_ad_{\widetilde{a}}-d_an_{\widetilde{a}}=-\pmb{(}n_{\widetilde{a}}(d_a-d_{\widetilde{a}})-d_{\widetilde{a}}(n_a-n_{\widetilde{a}})\pmb{)}$, \phantom{(1) }by the Cauchy-Schwarz inequality, 
$\textstyle |n_ad_{\widetilde{a}}-d_an_{\widetilde{a}}|\!\le\!\sqrt{S_{\widetilde{a}}}\;\!\eta$. So on $i\mR$, 
 \phantom{(1) }$
\textstyle 
{\scaleobj{0.81}{\cfrac{|n_ad_{\widetilde{a}}-d_an_{\widetilde{a}}|}{\sqrt{S_a}\sqrt{S_{\widetilde{a}}}}}}
\!\le\!{\scaleobj{0.81}{\cfrac{\eta}{\sqrt{S_a}}}}\le{\scaleobj{0.81}{\cfrac{\Lambda_a|a-\widetilde{a}|}{\sqrt{\delta_a}}}}.
$ 
\end{proof}

\noindent  Theorem~\ref{thm_quant} shows that $a\mapsto p_a$ is locally Lipschitz into $(\mS(H^\infty),d)$. 

\noindent {\bf The chordal distance in closed form.} 
Since $n_a,d_a$ are real on $(0,\infty)$, the reflection principle gives $n_a(\overline{s})\!=\!\overline{n_a(s)}$ and $d_a(\overline{s})\!=\!\overline{d_a(s)}$. So the chordal distance (the essential supremum over $i\mR$) is the supremum over $s\!=\!iy$, $y\!>\!0$ 
(the function is continuous on $i\mR\setminus\{0\}$, and extends continuously to $0$). For $y\!>\!0$, set $z\!:=\!\sqrt{iy}\!=\!\sqrt{y}\,\zeta$, where $\zeta\!:=\!e^{i\frac{\pi}{4}}$. Then $|z|^2\!=\!y$. We have 
$$
\textstyle
\begin{array}{rcl}
n_a d_{\widetilde{a}}-d_an_{\widetilde{a}}
\!\!\!&=&\!\!\!
{\scaleobj{0.81}{\cfrac{z\sinh(az)\sinh(\widetilde{a}z)\pmb{(}\cosh(az)-\cosh(\widetilde{a}z)\pmb{)}}{(1+z)^2\sinh z\cosh(az)\cosh(\widetilde{a}z)}}},\\[0.5cm]
|n_a|^2+|d_a|^2\!\!\!&=&\!\!\!{\scaleobj{0.81}{\cfrac{|\sinh(az)|^2\pmb{(}|\cosh(az)|^2+y|\sinh z|^2\pmb{)}}{|1+z|^2|\sinh z|^2|\cosh(az)|^2}}},
\end{array}
$$
and $|n_{\widetilde{a}}|^2+|d_{\widetilde{a}}|^2$ is given by the same formula with $\widetilde{a}$ in place of $a$.
 The factors $|1+z|$, $|\sinh(az)|$, $|\sinh(\widetilde{a}z)|$, $|\cosh(az)|$, $|\cosh(\widetilde{a}z)|$ disappear when we take the quotient, and we obtain 
$$
\textstyle 
\kappa(p_a,p_{\widetilde{a}})=\sup\limits_{y>0}\kappa(y),
$$
where with $z\!=\!\sqrt{y}\zeta$, 
$$
\textstyle 
\kappa(y):={\scaleobj{0.81}{\cfrac{\sqrt{y} |\sinh z|  |\cosh(az)-\cosh(\widetilde{a}z)|}{\sqrt{|\cosh(az)|^2+y|\sinh z|^2}\sqrt{|\cosh(\widetilde{a}z)|^2+y|\sinh z|^2}}}}.
$$

\noindent {\bf Determination of the location of the supremum.} 
Without loss of generality, let $a\!<\!\widetilde{a}$. For $y\!>\!0$, as above, $z\!=\!\sqrt{iy}\!=\!\sqrt{y}\,\zeta$. We have $|z|\!=\!\sqrt{y}$ and $x\!:=\!\text{Re}\;\!z\!=\!\sqrt{\frac{y}{2}}$. 
The function $\kappa$ is continuous on $(0,\infty)$, and $\kappa(y)\to0$ as $y\searrow 0$, since the numerator tends to $0$ and the denominator tends to $1$. 
The denominator of $\kappa(y)$ is at least $\sqrt{y|\sinh z|^2}\sqrt{y|\sinh z|^2}\!=\!y|\sinh z|^2$. So
$$
\textstyle
\kappa(y)\le{\scaleobj{0.81}{\cfrac{|\cosh(az)-\cosh(\widetilde{a}z)|}{|z||\sinh z|}}}.
$$
For $w\!=\!u+iv\in\mC$, where $u,v\in \mR$, we have 
$$
\textstyle 
\begin{array}{rcl}
|\cosh w|^2\!\!\!\!&=&\!\!\!\! (\sinh u)^2+(\cos v)^2\!=\!(\cosh u)^2-(\sin v)^2 \!\le\! (\cosh u)^2,\text{ and}\\[0.1cm] 
|\sinh w|^2\!\!\!\!&=&\!\!\!\!(\sinh u)^2+(\sin v)^2\!\ge\! (\sinh u)^2.
\end{array}
$$ 
As $\cosh$ is increasing on $[0,\infty)$, and since $a\!<\!\widetilde{a}$, we obtain  
$$
\textstyle 
\begin{array}{rcl}
|\cosh(az)-\cosh(\widetilde{a}z)|\!\!\!\!&\le&\!\!\!\!\cosh(ax)+\cosh(\widetilde{a}x)\le 2\cosh(\widetilde{a}x),
\text{ and}\\[0.1cm]
 |z||\sinh z|\!\!\!\!&\ge&\!\!\!\!\sqrt{y}\sinh x.
 \end{array}
$$
Thus
$$
\textstyle
\kappa(y)\le {\scaleobj{0.81}{\cfrac{2\cosh(\widetilde{a}x)}{\sqrt{y}\sinh x}}}.
$$
We claim that the right-hand side is a decreasing function of $y$. Firstly, $y\mapsto\frac{1}{\sqrt{y}}$ is decreasing. Also, $(0,\infty)\owns x\mapsto\frac{\cosh(\widetilde{a}x)}{\sinh x}$ is decreasing, 
as its derivative has the sign of $\widetilde{a}\sinh(\widetilde{a}x)\sinh x-\cosh(\widetilde{a}x)\cosh x<0$ (because $\widetilde{a}\!<\!1$), while $x\!=\!\sqrt{\frac{y}{2}}$ is increasing in $y$. 
Moreover, the right-hand side tends to $0$ as $y\to\infty$ (since $\widetilde{a}<1$).

Thus $\kappa$ is continuous on $(0,\infty)$, positive somewhere (as $a\!\neq\!\widetilde{a}$), and tends to $0$ both as $y\searrow0$ and as $y\to\infty$. 
So the supremum is a maximum, attained at some finite $y$. Furthermore, in any given example, 
once the right-hand side above has dropped below a value of $\kappa$ already found, the maximum can only be attained in the bounded interval of $y$'s up to that point, so one only needs to search a bounded interval.

\noindent {\bf Numerical example.} Let $a\!=\!\frac12$, $\widetilde{a}\!=\!\frac34$. At $y\!=\!100$ the bound above is $\approx\!0.034$, and it is decreasing in $y$, so $\kappa(y)\!<\!0.035$ for all $y\!\ge\!100$. 
Since $\kappa$ takes values of about $0.13$ at moderate $y$ (see below), the maximum is attained in $(0,100]$.  
A plot of $\kappa$ is shown below. It was determined numerically that $\kappa$ attains its maximum at $y_*\!\approx\!3.9302$, with value $\kappa(y_*)\!\approx\!0.13501$.
\begin{figure}[H]
     \center
     \includegraphics[width=8.1 cm]{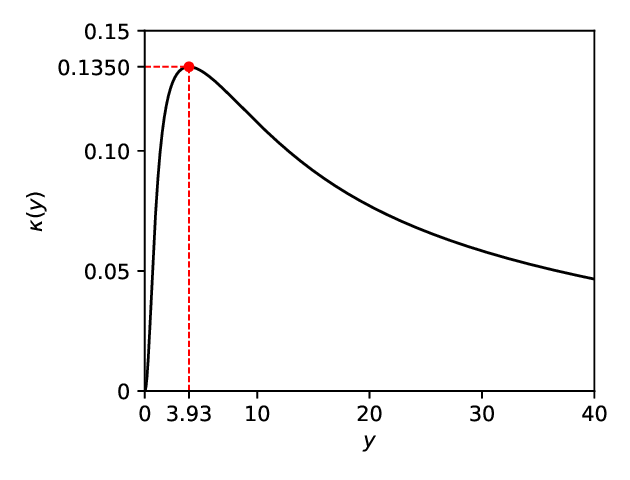}
     
     {\small Plot of the chordal distance $y\mapsto \kappa(y)$ for $a\!=\!\frac12$ and $\widetilde{a}\!=\!\frac34$.\\ The maximum $\approx\!0.1350$ is attained at $y_*\!\approx\!3.93$.\phantom{The maxi}}
\end{figure}

\noindent Moreover, with $f\!:=\!n_a^*n_{\widetilde{a}}+d_a^*d_{\widetilde{a}}$ on $i\mR$, the function $\text{Re}\;\!f$ is continuous on $i\mR$, and  
$\text{Re}\;\!f(iy)\to a\widetilde{a}\!=\!\frac{1}{2}\cdot \frac{3}{4}\!=\!\frac38$ as $y\to0$ (since $n_a(0)\!=\!a$, $n_{\widetilde{a}}(0)\!=\!\widetilde{a}$ and $d_a(0)\!=\!d_{\widetilde{a}}(0)\!=\!0$), and $\text{Re}\;\!f(iy)\to 1$ as $y\to\infty$ (since $n_a,n_{\widetilde{a}}\to0$ and $d_a,d_{\widetilde{a}}\to1$ as $|s|\to\infty$).
Its graph (plot below) shows that $\text{Re}\;\!f(iy)\!>\!0$ for all $y$. Its minimum is  $\approx\!0.187$, attained at $y_{\star}\!\approx\!0.51$.  For large $y$ it stays close to $1$. Hence $\text{Re}\;\!f\!\ge\!\epsilon$ on $i\mR$ for $\epsilon\!=\!0.18\!>\!0$. By Lemma~\ref{lem_re}, $f$ is invertible and the index condition $\iota(f)\!=\!0$ is satisfied. 
Note that this example, where $|a-\widetilde{a}|\!=\!\frac14$, lies outside the range of Theorem~\ref{thm_quant}. 
Indeed, by Proposition~\ref{28_1_25_1011}, $\delta_a\!\le\!|n_a(s)|^2+|d_a(s)|^2$ for every $s\in\overline{\mC_+}$, and evaluating at $s\!=\!\frac{i}{2}$ gives $\delta_a\!\le\!0.125$ for $a\!=\!\frac12$.
 With 
 $$
 \textstyle \Lambda_a\!=\sqrt{4+\frac{8}{a^2}}\!=\!6,
 $$
  we have 
  $$
  \textstyle
  {\scaleobj{0.81}{\cfrac{\sqrt{\delta_a}}{2\Lambda_a}}}\le 0.03,
  $$ 
much smaller than $|a-\widetilde{a}|\!=\!\frac14\!=\!0.25$. Thus, instead of appealing to Theorem~\ref{thm_quant}, we verify the index condition directly, using the numerically computed graph of $\text{Re}\;\!f$. Consequently, 
$$
\textstyle
d(p_{\frac{1}{2}},p_{\frac{3}{4}})\!=\!\kappa(p_{\frac{1}{2}},p_{\frac{3}{4}})\!\approx\! 0.1350.
$$
\begin{figure}[H]
     \center
     \includegraphics[width=12 cm]{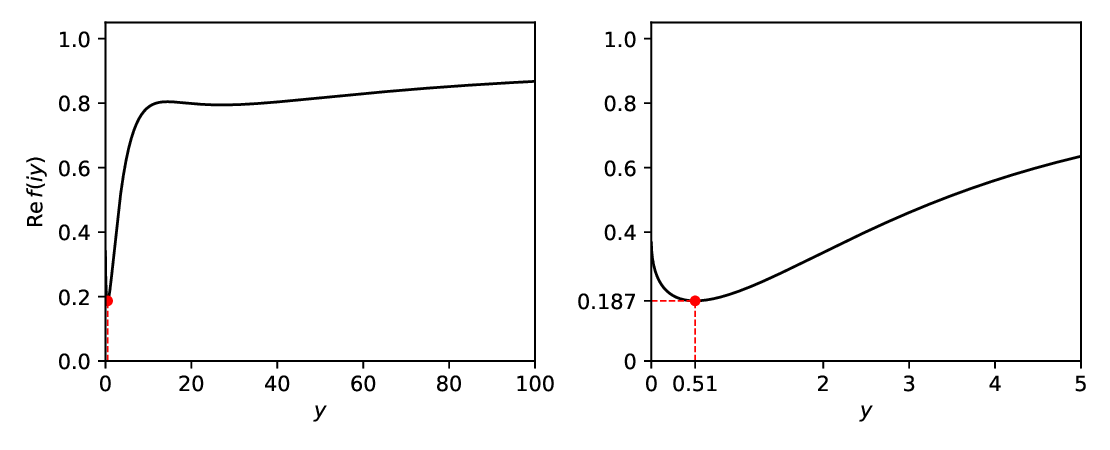}
     
     {\small Plot of $y\mapsto\text{Re}\;\!f(iy)$ for $a\!=\!\frac12$, $\widetilde{a}\!=\!\frac34$, for $y\in[0,100]$ (left), 
     and \\a close-up for $y\in[0,5]$ (right). 
     The minimum $\approx\!0.187$ occurs at $y\!\approx\!0.51$.}
\end{figure}

\noindent For other $\widetilde{a}$ (with $a\!=\!\frac{1}{2}$) the same computation gives the following (in each case $\text{Re}\;\!f$ has a positive minimum, so $\iota\!=\!0$ and $d\!=\!\kappa$):
\begin{table}[H]
\centering
{\footnotesize
\setlength{\tabcolsep}{4pt}
\renewcommand{\arraystretch}{1.15}
\begin{tabular}{|c|c|c|c|c|c|c|}
\hline
\rowcolor{gray!25}
$\widetilde{a}$ & $0.3$ & $0.4$ & $0.5$ & $0.6$ & $0.7$ & $0.8$ \\ \hline
$d(p_{\frac{1}{2}},p_{\widetilde{a}})\rule[-1.5ex]{0pt}{4.4ex}$ & $0.0698$ & $0.0393$ & $0$ & $0.0479$ & $0.1040$ & $0.1678$ \\ \hline
$\min\limits_{y}\text{Re}\;\!f(iy)\rule[-1.8ex]{0pt}{3.8ex}$ & $0.075$ & $0.100$ & $0.125$ & $0.150$ & $0.174$ & $0.199$ \\ \hline
\end{tabular}}
\end{table}

 \noindent {\bf Acknowledgements:} The author thanks Raymond Mortini (Universit\'e du Luxembourg) 
 for useful discussions, and the two anonymous referees for their valuable comments.

\smallskip 

\noindent {\bf Conflict of interests:} There are no competing interests to declare.

\end{document}